\documentclass[12pt,a4paper]{amsart}

\input{prepictex}
\input{pictex}
\input{postpictex}
\theoremstyle{plain}

\newtheorem*{thrm*}{Theorem}

\newtheorem*{prop*}{Proposition}

\newtheorem*{lemma*}{Lemma}

\newtheorem*{cor*}{Corollary}
\theoremstyle{remark}

\theoremstyle{remark}
\newtheorem*{example}{Example}

\newtheorem*{examples}{Examples}

\newcommand{\R}{\mathbb R}                                   

\newcommand{\set}[1]{\left\{#1\right\}}                      
\newcommand{\abs}[1]{\left\lvert#1\right\rvert}              

\DeclareMathOperator{\Cl}{Cl}                                
\DeclareMathOperator{\Lk}{Lk}                                
\DeclareMathOperator{\conv}{conv}                            
\DeclareMathOperator{\dimm}{dim}                             
\DeclareMathOperator{\St}{St}                                
\newcommand{\simp}[1]{$#1$\nobreakdash-simplex}              
\newcommand{\simps}[1]{$#1$\nobreakdash-simplices}           
\newcommand{\dimn}[1]{$#1$\nobreakdash-dimensional}          
\newcommand{\manif}[1]{$#1$\nobreakdash-manifold}            
\newcommand{\sph}[1]{$#1$\nobreakdash-sphere}                

\theoremstyle{remark}

\theoremstyle{plain}
\newtheorem*{lemma1}{Lemma 1}
\newtheorem*{lemma2}{Lemma 2}
\newtheorem*{defn}{Definition}

\begin{document}

\title{Euler Characteristic in Odd Dimensions}
\author{Colin MacLaurin and Guyan Robertson}
\address{Department of Mathematics  \\
        University of Newcastle\\  NSW  2308\\ AUSTRALIA}
\email{colin.maclaurin@studentmail.newcastle.edu.au}
\email{guyan.robertson@newcastle.edu.au }

\begin{abstract}
It is well known that the Euler characteristic of an odd dimensional compact manifold is zero.
An Euler complex is a combinatorial analogue of a compact manifold.  
We present here an elementary proof of the corresponding result for Euler complexes.
\end{abstract}

\maketitle

\section{Background}

We begin by recalling background information on simplicial complexes
and the motivating result from topology. Only finite simplicial complexes will be considered, and
we include this condition in the definitions.

\subsection*{Geometric Simplicial Complexes}
Let $S$ be a finite subset of $\R^n$. An \emph{affine combination} of the points $p_i\in S$ is a point $x=\sum\lambda_ip_i$, where $\lambda_i\in\R$, $\sum \lambda_i=1$. The set $S$ is called \emph{affinely independent} if no point in $S$ is an affine combination of the other points. A \emph{convex combination} is an affine combination with non-negative coefficients $\lambda_i$. The \emph{convex hull}, $\conv S$ is the set of all convex combinations of the points in $S$.

A \emph{(geometric) \simp{k}} is the convex hull of a set of $k+1$ affinely independent points
in $\R^n$. We say $k$ is the \emph{dimension} of the simplex. In $\R^n$ there are simplices ranging from dimension $-1$ through to $n$.

Let $S\subset\R^n$ be finite and affinely independent (so $\sigma=\conv S$ is a simplex). If $T\subset S$ then $T$ is also affinely independent, so the convex hull $\tau=\conv T$ is a simplex as well. We say $\tau$ is a \emph{face} of $\sigma$.

\begin{example}
A \simp{2} in $\R^2$ is a triangle, and has eight faces: the triangle itself, its three edges, three vertices, and also the empty set $\emptyset$.
\end{example}

A \emph{(geometric) simplicial complex} is a finite collection $C$ of simplices such that
\begin{enumerate}
\item [(i)]  if $\sigma\in C$, and $\tau$ is a face of $\sigma$ then $\tau\in C$\,;
\item [(ii)] if $\sigma,\tau\in C$ then  $\sigma\cap\tau$ is a face of both $\sigma$ and $\tau$.
\end{enumerate}

Note that $\sigma$ and $\tau$ may be be disjoint in (ii), since $\emptyset$ is a face of every simplex.
The dimension of $C$ is the highest dimension of a simplex contained in $C$.

Let $C$ be an \dimn{n} simplicial complex. The \emph{polyhedron} or \emph{underlying space} of $C$, written $\abs{C}$, is the topological space formed by taking the subset of $\R^n$ consisting of points which are covered by simplices in $C$,  with the usual topology inherited from $\R^n$.

\subsection*{Abstract Simplicial Complexes}
An \emph{abstract simplicial complex} is a finite collection $X$ of finite sets 
with the following property :
\begin{itemize}
\item if $\alpha\in X$ and $\beta\subset\alpha$ then $\beta\in X$.
\end{itemize}

A set $\alpha\in X$ is called an \emph{(abstract) simplex} of dimension $\dimm\alpha=\#\alpha-1$.
If $\beta\subset\alpha$, we say that $\beta$ is a face of $\alpha$.
If $k\ge 0$ then $X^k$ will denote the set of $k$-dimensional simplices of $X$. 
Elements of $X^0$ are called vertices and the dimension of $X$ is the maximum dimension of a simplex in $X$. A simplicial complex $X$ is said to be \emph{pure} if all maximal simplices of $X$ have the same dimension $n$. 

The definition of an abstract simplicial complex is elegant and useful for combinatorial calculations. However for the sake of intuition, it is convenient to represent it geometrically.
A \emph{geometric realisation} of an abstract simplicial complex $X$ is a geometric simplicial complex $C$ with a bijection $\phi$ from the vertices of $X$ to the vertices of $C$ such that $\alpha\in X$ if and only if $\conv\phi(\alpha)\in C$.
It is well known that a \dimn{k} abstract simplicial complex can always be realised in $\R^{2k+1}$.

\subsection*{Topology}

An \emph{\manif{n}} is a topological space $M$ in which every point of $M$ has a neighbourhood homeomorphic to an open $n$-ball in $R^n$.

\begin{examples}
For $n\ge 1$,  $\R^n$ itself is an \manif{n}, as is the \sph{n} $S^n$.
\end{examples}

\subsection*{Euler Characteristic, Triangulation}

Let $X$ be an \dimn{n}  simplicial complex, and let $s_k$ be the number of \dimn{k} simplices of $X$. The \emph{Euler characteristic} is defined by
\[ \chi(X)=\sum_{k=0}^n (-1)^k s_k. \]
This number is also the Euler characteristic of any geometric realisation of $X$.

\begin{example}
Let $T$ be the \dimn{2} simplicial complex which is the \emph{surface} of a tetrahedron. $T$ contains $4$ vertices, $6$ edges, and $4$ faces. Hence $\chi(T)=4-6+4=2$. \end{example}

A geometric simplicial complex $C$ is a \emph{triangulation} of a topological space $M$ if $\abs{C}$ is homeomorphic to $M$. We say $M$ is \emph{triangulable} if it has a triangulation. Intuitively, a triangulation gives a topological space a combinatorial structure.
It is known that if $K$ and $L$ are triangulations of the same topological space $M$ then $\chi(K)=\chi(L)$.
If $M$ is triangulable, the Euler characteristic of $M$ may therefore be defined to be the Euler characteristic of any triangulation of $M$. 
 
\begin{example}
The surface $M$ of an \simp{(n+1)} is a triangulation of the \sph{n} $S^n$. It is not hard to see that $\chi(M)=1+(-1)^n$, and hence $\chi(S^n)=1+(-1)^n$.
\end{example}

It is known that all compact manifolds of dimensions 2 and 3 can be triangulated, but that there exist compact manifolds of dimension 4 which cannot be triangulated \cite{BL}. Now the Euler characteristic classifies the compact orientable manifolds of dimension 2.  However this fails dramatically in dimension 3. In fact every 3-dimensional compact manifold has Euler characteristic zero. More generally, the Euler characteristic can be defined by other methods for any compact manifold and one obtains:

\begin{thrm*}
The Euler characteristic of an odd-dimensional compact manifold is zero.
\end{thrm*}

There is a proof of this fact in \cite[Cor 3.37]{Hat}. It is nontrivial and uses Poincar\'{e} duality. The book \cite{Hat} is an excellent reference for the topological background.
The combinatorial aspects are also explored in further detail in \cite{Mun}.

\section{Euler Complexes}

We consider a class of simplicial complexes, known as Euler complexes, which 
embody the defining property of a manifold.
Our result may be well known, but we know of no reference to it.  It gives an intuitive insight into why the manifold result is true. Since topology plays no role, we will use the language of abstract simplicial complexes.

Let $\sigma$ be a simplex in a simplicial complex $X$. There are two important subcomplexes of $X$ associated to $\sigma$. The \emph{star} of $\sigma$ is the set of simplices containing $\sigma$, and the \emph{link} of $\sigma$ consists of all faces of simplices in the star that do not intersect $\sigma$. We use the following notation for these concepts:

\begin{align*}
\St\sigma &= \set{\tau\in X\, \colon\sigma\subseteq\tau} \\
\Lk\sigma &= \set{\tau\in\Cl(\St\sigma)\,\colon\sigma\cap\tau=\emptyset}.
\end{align*}

Here $\Cl$ denotes the \emph{closure} of $\St\sigma$. The closure of a subcomplex $Y$ of a simplicial complex consists of all faces of simplices in $Y$.

\begin{defn} \cite{Sat}
A pure simplicial complex is an \emph{Euler complex} if the Euler characteristic of the link of every simplex is equal to the Euler characteristic of the sphere of the same dimension.
\end{defn}

\begin{figure}[htbp]
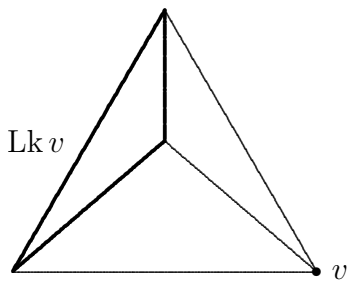
\label{fig1}  
{}\hfill
\beginpicture
\setcoordinatesystem units <1cm,1.732cm>  
\setplotarea  x from -2.5 to 2.5,  y from -1.5 to 1.5
\setlinear \plot 2 -1  -2 -1  0 1  2 -1  /
\plot 0 0  2 -1 /
\plot 0 0  -2 -1 /
\plot 0 0  0 1 /
\put {$\Lk v$} at -1.7 0
\put {$v$} at 2.3 -1
\put {$_\bullet$} at 2 -1
\setplotsymbol({$\cdot$})
 \plot -2 -1  0 1  0 0  -2 -1 /  
\endpicture
\hfill{}
\caption{A vertex and its link on the surface of a tetrahedron.}
\end{figure}

In Figure \ref{fig1}, $X$ is the surface of a tetrahedron. This is an Euler complex. The link of a vertex $v$ is a triangle. The link of an edge is the set of two vertices not lying on the edge.

\begin{thrm*}
The Euler characteristic of an Euler complex of odd dimension is zero.
\end{thrm*}

From now on, let $X$ be an Euler complex of dimension $n$, and introduce the following notation:
\begin{itemize}
\item $\check s_l(Y)$ is the number of \simps{l} contained in a subcomplex $Y\subset X$,
\item $\hat s_l(\sigma)$ is the number of \simps{l} in $X$ which contain a simplex $\sigma\in X$.
\end{itemize}

\begin{lemma1}
Let $0\le k \le l \le n$. Then
\begin{equation}\label{one}
\sum_{\sigma\in X^k}\hat s_l(\sigma)=\binom{l+1}{k+1}\check s_l(X).
\end{equation}
\end{lemma1}

\begin{proof} 
If $(\sigma,\tau)\in X^k\times X^l$ define
$$\delta_{\sigma,\tau}=\begin{cases}
1& \text{if $\sigma\subseteq \tau$},\\
0& \text{otherwise.}
\end{cases}$$
Now the number of \simps{k} $\sigma$ contained in a fixed \simp{l} $\tau$ is $\binom{l+1}{k+1}$.
Therefore
\begin{equation*}
\sum_{\sigma\in X^k}\hat s_l(\sigma)=\sum_{\sigma\in X^k}\sum_{\tau\in X^l}\delta_{\sigma,\tau} 
=\sum_{\tau\in X^l}\sum_{\sigma\in X^k}\delta_{\sigma,\tau}  
=\sum_{\tau\in X^l}\binom{l+1}{k+1} 
=\binom{l+1}{k+1}\check s_l(X)\,.
\end{equation*}
\end{proof}

\begin{lemma2}
Let $\sigma$ be a \simp{k}. The number of \simps{l} containing $\sigma$ is equal to the number of \simps{(l-k-1)} in $\Lk\sigma$; that is,
\begin{equation}\label{two}
 \hat s_l(\sigma)=\check s_{l-k-1}(\Lk\sigma). 
\end{equation}
\end{lemma2}

\begin{proof}
Let $\tau$ be a \simp{(l-k-1)} in $\Lk\sigma$.
Then, by definition, $\tau$ is a face of some simplex $\rho$ with $\rho \supset \sigma$ and  
$\sigma\cap\tau=\emptyset$.
This means that $\sigma \cup \tau$ is a face of $\rho$ and contains $l-k+k+1=l+1$ vertices.
Therefore $\sigma \cup \tau$ is a \simp{l}.
Now the map $\tau \mapsto \sigma\cup\tau$ is a bijection from the set of \simps{(l-k-1)} in $\Lk\sigma$ onto the set of \simps{l} containing $\sigma$, with inverse map
$\rho \mapsto \rho \setminus \sigma$.
\end{proof}

\begin{proof}[Proof of the Theorem] If $\sigma\in X^k$ then
\begin{align*}
\chi(\Lk\sigma) &= \sum_{p=0}^{n-k-1}(-1)^p\check s_p(\Lk\sigma) \\
&= \sum_{l=k+1}^n(-1)^{l-k-1}\check s_{l-k-1}(\Lk\sigma) & (\text{writing } l=p+k+1) \\
&= \sum_{l=k+1}^n(-1)^{l-k-1}\hat s_l(\sigma) & \text{(by (\ref{two}))}\,. \\
\end{align*}
Now, since $X$ is an Euler complex, $\chi(\Lk\sigma)=1+(-1)^k$. Therefore summing over
all \simps{k} $\sigma$ gives
\begin{align*}
(1+(-1)^k)\check s_k(X) &= \sum_{\sigma\in X^k}\sum_{l=k+1}^n(-1)^{l-k-1}\hat s_l(\sigma) \\
&= \sum_{l=k+1}^n(-1)^{l-k-1}\sum_{\sigma\in X^k}\hat s_l(\sigma)  \\
&= \sum_{l=k+1}^n(-1)^{l-k-1}\binom{l+1}{k+1}\check s_l(X)\,, & \text{(by (\ref{one}))\,.} \\
\end{align*}
Now take the sum from $k=0$ to $n-1$ of both sides. Since $n$ is odd, the left hand side becomes $2\sum_{k\text{ even}}\check s_k(X)$. Therefore

\begin{equation}\label{main}
2\sum_{k\text{ even}}\check s_k(X) = \sum_{k=0}^{n-1} \left(\sum_{l=k+1}^n (-1)^{l-k-1}\binom{l+1}{k+1}\check s_l(X)\right). 
\end{equation}
\\
The coefficient of $\check s_l(X)$ on the right hand side is 
\begin{align*}
& \sum_{k=0}^{l-1}(-1)^{l-k+1}\binom{l+1}{k+1} \\
=& (-1)^{l+1}-(-1+1)^{l+1}+1 & \text{(by the Binomial Theorem)} \\
=& \begin{cases}
0 & l \text{ even} \\
2 & l \text{ odd}.
\end{cases}
\end{align*}
Therefore equation (\ref{main}) becomes
\[ 2\sum_{k\text{ even}}\check s_k(X)=2\sum_{l\text{ odd}}\check s_l(X), \]
from which we obtain $\chi(X)=0$.
\end{proof}


\begin{thebibliography}{MUN} 
\bibitem [BL]{BL} A. Bj\"orner and F. Lutz, Simplicial manifolds, bistellar flips and a $16$-vertex triangulation of the Poincar\'e homology $3$-sphere, Experiment. Math. {\bf 9}
(2000), 275--289.
\bibitem [Hat]{Hat} A. Hatcher, \emph{Algebraic Topology}, Cambridge University Press,
Cambridge 2002.
\bibitem [Mun]{Mun} J.R. Munkres, \emph{Topology: A First Course}, Prentice-Hall, Englewood Cliffs, N.J., 1975.
\bibitem [Sat]{Sat} H. Sato, \emph{Algebraic Topology: An Intuitive Approach}, American Mathematical Society, Providence, R.I., 1999.
\end{thebibliography}
\end{document}